\documentclass[reqno]{amsart}
\usepackage{mathtools,amssymb,wasysym}
\usepackage{mathrsfs}
\usepackage{enumitem}
\usepackage[usenames, dvipsnames]{color}

\usepackage[numbers,sort&compress]{natbib}
\usepackage{hyperref}
\usepackage{esint}
\usepackage{verbatim}

\theoremstyle{plain}
\newtheorem{theorem}{Theorem}[section]
\newtheorem{lemma}[theorem]{Lemma}
\newtheorem{corollary}[theorem]{Corollary}

\theoremstyle{definition}
\newtheorem{definition}[theorem]{Definition}

\theoremstyle{remark}
\newtheorem{remark}[theorem]{Remark}
\newtheorem*{rem*}{Remark}

 \makeatletter
 \newcommand{\norm}{\@ifstar{\@normb}{\@normi}}
 \newcommand{\@normb}[2]{\left\Vert{#1}\right\Vert_{#2}}
 \newcommand{\@normi}[2]{\Vert{#1}\Vert_{#2}}
 \makeatother

 \RequirePackage{bm}

 \DeclareMathOperator{\supp}{supp}

\newcommand{\liu}{(-\Delta)^{-1}|u|^2}
\newcommand{\R}{\mathbb{R}}

\newcommand{\M}{\mathcal{M}}
\newcommand{\MP}{\mathcal{M}_{\Psi}}
\newcommand{\li}{(-\Delta)^{-1}}
\newcommand{\hu}{\widehat{u_{+}}}
\newcommand{\p}{\partial}
\newcommand{\asy}{\widehat{u_{+}}}
\newcommand{\F}{\mathcal{F}}
\newcommand{\U}{\mathcal{U}}

\makeatletter
\def\@tocline#1#2#3#4#5#6#7{\relax
  \ifnum #1>\c@tocdepth 
  \else
    \par \addpenalty\@secpenalty\addvspace{#2}%
    \begingroup \hyphenpenalty\@M
    \@ifempty{#4}{%
      \@tempdima\csname r@tocindent\number#1\endcsname\relax
    }{%
      \@tempdima#4\relax
    }%
    \parindent\z@ \leftskip#3\relax \advance\leftskip\@tempdima\relax
    \rightskip\@pnumwidth plus4em \parfillskip-\@pnumwidth
    #5\leavevmode\hskip-\@tempdima
      \ifcase #1
       \or\or \hskip 1em \or \hskip 2em \else \hskip 3em \fi%
      #6\nobreak\relax
    \hfill\hbox to\@pnumwidth{\@tocpagenum{#7}}\par
    \nobreak
    \endgroup
  \fi}
\makeatother

\usepackage{tcolorbox}

\allowdisplaybreaks

\begin{document}

\title[Hartree with Coulomb]{Modified Wave operators for the Hartree equation with repulsive Coulomb potential}

\author[W. Huang]{Wenrui Huang}
\address{Department of Mathematics, Brown University, 151 Thayer Street, Providence, RI 02912, USA}
\email{wenrui\_huang@brown.edu }


\begin{abstract}
We study the final state problem for the Hartree equation with repulsive Coulomb potential:
\[i\partial_t u+\frac{1}{2}\Delta u-\frac{1}{|x|}u=((-\Delta)^{-1}|u|)^2u\]
We show the work in \cite{KaMi} can be extended to the Hartree nonlinearity: Given a prescribed asymptotic profile, we construct a unique global solution scattering to the profile. In particular, the existence of the modified wave operators is obtained for sufficiently localized small scattering data.
\end{abstract}

\maketitle 
\section{Introduction}
In this paper, we consider the final state problem of the Hartree equation with repulsive Coulomb potential in $\R^3$:
\begin{equation}\label{HCoulomb}
	\begin{split}
		i\partial_t u-Hu=\left(\frac{1}{|x|}*|u|^2\right)u ,\qquad x\in\R^3,\qquad t\in \R
	\end{split}
\end{equation}
where
\[ Hu=-\frac{1}{2}\Delta u+\frac{K}{|x|}u, \qquad H_0u:=-\frac{1}{2}\Delta u \]

The study of the operator $H=-\Delta+\frac{K}{|x|}$ with Coulomb potential originates from both physical and mathematical interests. In particular, when $K<0$, This operator provides a quantum description of the Coulomb force between two charged particles and corresponds to attractive long-range potential due to a positively charged atomic nucleus.

The mathematical interest in these equations comes from the operator theory with a long-range potential and the solution's dispersive behavior. Observing that $|x|^{-1}\in L^{2}({\R^3})+L^{\infty}(\R^3)$, we know that $H$ is essentially self-adjoint on $C_0^{\infty}(R^3)$ and self-adjoint. Moreover, when $K>0$, we refer to \cite{black2024pointwisedecayradialsolutions} for recent work studying dispersive estimates with radial data using distorted Fourier transform. Also, the nonlinear equation \eqref{HCoulomb} has been studied extensively in the literature. We refer to \cite{ChGlassey}\cite{Dias}\cite{Hayashi:1987aa}  for the global existence and a decay rate for the solution when $K\geq 0$. 

The natural question is: what is the asymptotic behavior?  However, there is a paucity of literature on the asymptotic behavior of $\eqref{HCoulomb}$. Without the Coulomb potential, Ginibre and Ozawa \cite{GO93} proved the modified scattering result for the Hartree equation. Nakanishi \cite{N02A} and \cite{N02B} proved the existence of modified wave operators for the Hartree equation.

In this paper, we consider the existence of modified wave operators of the Hartree equation with the repulsive Coulomb potential. We hope our final state result can be helpful for the study of the modified scattering problem of \eqref{HCoulomb}.

From now on, we will assume $K=1$. Our main result is the following: Given a scattering datum $u_{+}$, we construct a unique global solution scattering to a prescribed asymptotic profile $u_p(t)$ defined in \eqref{asymprfile}.

Our main result is the following:

\begin{theorem}\label{MainTheorem}
      Suppose $\frac{1}{4}<b<\frac{1}{2}$, $\frac{2}{q}+\frac{3}{r}=\frac{3}{2}$ , $\widehat{u_+}\in H^1(\R^3)$ and $0\notin \mathrm{supp}\; \widehat{u_{+}}$ and $\|\li|\asy|^2 \|_{L^\infty}$ be small enough. Then there exists a unique solution $u\in C(\R;L^2(\R^3))$  to (1) satisfying, for any $(q,r)$,
      $$ \|u(t)-u_p(t)\|_{L^2(\R^3)}+\|u-u_p\|_{L^q([t,\infty);L^r(\R^3)}\lesssim t^{-b}, \quad t\rightarrow +\infty
      $$
  \end{theorem} 
\begin{remark}\label{SoblevEmbed} From Young's inequality and Sobolev embedding, we have 
\[ \|(-\Delta)^{-1}|f|^2\|_{L^{\infty}(\R^3)}\lesssim \| f\|^2_{L^2(\R^3)\cap L^4(\R^3)}\lesssim \| f\|^2_{H^1(\R^3)}
\]
	
\end{remark}  
\begin{remark}
	The condition $0\notin \supp \asy$ is due to Lemma \ref{ControlU1U2} . It is interesting to ask if we can remove or replace this condition with some weight.
\end{remark}

As an immediate consequence, we obtain the existence of the modified wave operator:
  
 \begin{corollary}[{Modified wave operator}]
    There exists $\varepsilon_0>0$ such that for any $0<\varepsilon\leq \varepsilon_0$, we have the existence of the modified wave operator
    $$ W^{+}_{\Psi}:u_+ \mapsto u(0)
    $$
    which is defined from $\mathcal{F}H^1\cap \{u\in L^2(\R^3)| \,0\notin \mathrm{supp}\;  \hat{u} \; \mathrm{and}\; \|(-\Delta)^{-1}|\hat{u}|^2\|_{L^\infty}\leq \varepsilon \}$ into $L^2(\R^3)$.
\end{corollary}

The work of \cite{KaMi} highly influences our proof of Theorem \ref{MainTheorem}. In \cite{KaMi}, the authors study the final state problem for the nonlinear Schrodinger equation with a critical long-range nonlinearity and a long-range potential. Their class of potential includes a repulsive long-range potential with a short-range perturbation. Theorem \ref{MainTheorem} also holds in this class. However, we only treat the case of  Coulomb potential since it is more interesting.

\textbf{Acknowledgement.} The author thanks Prof. Beno\^it Pausader for several comments and helpful discussions.

\section{preliminaries}
\subsection{Notations and definitions} We introduce some notations and definitions.
Let $\R^3$ denote the standard Euclidean space of points $x=(x^1,x^2,x^3)$ and let $\langle x\rangle=(1+x^2)^{1/2}$. We write $A\lesssim B$ if there exists a constant $C>0$ such that $A\lesssim CB$, where $C$ is independent of $A$ and $B$. We write $A\approx B$ if $A\lesssim B$ and $B\lesssim A$. We often suppress the space variables $\R^3$ when it involves norms, i.e.
\begin{equation*}
	  \|f\|_{L^r}=\|f \|_{L^r(\R^3)}, \qquad \| f\|_{H^s}=\|f \|_{H^s(\R^3)}  
\end{equation*}
And we use $\hat{f}(\xi)=\mathcal{F}f(\xi)$ to denote the Fourier transform.

\begin{definition}\label{DefinitionOperator}
Given smooth functions $f(x),\Psi(t,x)$ and $u_{+}(x)$, we define the following operators and functions.
	\begin{enumerate}
		\item $[D(t)f](x)=(it)^{-\frac{3}{2}}f(\frac{x}{t}) $ 
		
		\item $[\mathcal{M}(t){f}](x)=e^{\frac{i|x|^2}{2t}}f(x) $
		
		\item $[\MP(t)f](x)=e^{{i\Psi(t,x)}}f(x) $
		
		\item $W(t,x)=e^{[-i\li |\hu|^2](x)\log t}\hu(x)$
		
		\item $[F(u)](t,x):=[\liu](t,x)u(t,x)  $
		
	\end{enumerate}
	\end{definition}

Throughout this paper, we will frequently use these operators. Since $\M(t)$ and $\MP(t)$ are multiplication by a modulus one function, they are isometries on $L^p(\R^3)$.  And since $D(t)$ are dilations, we have 
\[ \norm{ D(t)f}{L^r}=|t|^{-3(\frac{1}{2}-\frac{1}{r})}\|f\|_{L^r}  \]

Under definition \eqref{DefinitionOperator}, the equation \eqref{HCoulomb} simply reads:
\begin{align*}
	i\p_{t}u-Hu=F(u)
\end{align*}

\begin{lemma}\label{Wequation} Under Definition \eqref{DefinitionOperator}, we have 
\begin{equation}\label{EquationW}
	i\partial_t W=\frac{F(W(t))}{t}
\end{equation} 
 \begin{equation}\label{EquationUp}
	F(u_p)(t,x)=\MP(t)D(t)t^{-1}F(W(t))=\MP(t)D(t)(i\partial_t W)
\end{equation}	
\end{lemma}

\begin{proof}
(1). Observe that $|W|=|\hu|$, for the first quality we compute directly
\begin{equation*}
	\begin{split}
		\partial_t W(t,x)=-\frac{i}{t}\li |\hu|^2W(t,x)=-\frac{i}{t}\li|W|^2W(t,x)
	\end{split}
\end{equation*}
(2). For the second equality, observe that 
\begin{equation*}
	\begin{split}
		&|u_p|^2= t^{-3} \left|\hu(\frac{x}{t})  \right|^2 \\
	\end{split}
\end{equation*}
Then, we can compute 
\begin{align*}
	\li|u_p|^2(t,x)&=t^{-3} \int_{\R^3} \frac{|\hu(\frac{y}{t})|^2}{|x-y|}dy \\
	&=\frac{1}{t}\int_{\R^3} \frac{|\hu(y)|^2}{|\frac{x}{t}-y |}dy \\
	&=\frac{1}{t}[\li|\hu|^2](\frac{x}{t})
\end{align*}
and 
\begin{equation*}
F(u_p)=\MP(t)\frac{1}{t}[\li|\hu|^2](\frac{x}{t}) D(t)W(t,x)=\MP(t)D(t)t^{-1}F(W(t))	
\end{equation*}
 
\end{proof}

\begin{lemma}[Nonlinear Estimates]\label{NonlinearEstimates}
\begin{equation*}
\begin{split}
		\|e^{-i(-\Delta)^{-1}|u|^2} u \|_{H^1}&\lesssim \|u\|_{H^1}(1+\|u\|^2_{H^1}) \\
\|(-\Delta)^{-1}|u|^2 e^{-i(-\Delta)^{-1}|u|^2}u \|_{H^{1}}&\lesssim \|u\|_{H^1}(1+\|u\|_{H^1}^2+\|u\|_{H^1}^4)
	\end{split}
\end{equation*}

\end{lemma}	
\begin{proof}
This lemma follows from Leibniz rule, Sobolev embedding, and Hardy-Littlewood-Sobolev inequality. Since the first line is more direct, let's consider the second line. By Remark \ref{SoblevEmbed}, we have
\begin{equation*}
	\begin{split}
	\|(-\Delta)^{-1}|u|^2 e^{-i(-\Delta)^{-1}|u|^2}u \|_{L^2}=\| (-\Delta)^{-1}|u|^2 u \|_{L^2}     
\leq \| u\|_{H^1}^2\|u\|_{L^2}
	\end{split}
\end{equation*}
Computing the derivative, we get:
\begin{equation*}
	\begin{split}
\nabla \{ (-\Delta)^{-1}|u|^2 e^{-i(-\Delta)^{-1}|u|^2}u \}	&=[\nabla (-\Delta)^{-1}|u|^2]e^{-i(-\Delta)^{-1}|u|^2}u \\
&+	(-\Delta)^{-1}|u|^2 [(-i)\nabla (-\Delta)^{-1}|u|^2] e^{-i(-\Delta)^{-1}|u|^2}u \\
&+ (-\Delta)^{-1}|u|^2]e^{-i(-\Delta)^{-1}|u|^2}\nabla u
	\end{split}
\end{equation*}
And by H\"older's inequality and Sobolev embedding, we have 
\begin{equation*}
	\begin{split}
		\|\nabla (-\Delta)^{-1}|u|^2]e^{-i(-\Delta)^{-1}|u|^2}u\|_{L^2}\leq \| \nabla (-\Delta)^{-1}|u|^2\|_{L^6}\|u\|_{L^3}\lesssim \| |u|^2 \|_{L^2}\| u\|_{L^4}\lesssim \|u\|_{H^1}^3
	\end{split}
\end{equation*}
The second and third terms follow directly. 
\end{proof}

\subsection{Strichartz Estimates} It is well known that the Strichartz estimate is very useful in studying nonlinear dispersive equations. First  a pair $(q,r)\in \R^2$ is said to be admissible if 
  \begin{equation}
  	2\leq q,r\leq \infty, \quad \frac{2}{q}+\frac{3}{r}=\frac{3}{2}
  \end{equation} 
  
In the repulsive Coulomb potential case, Mizutani \cite{Mizutani} recently proved the global-in-
time Strichartz estimate, where the proof employs several techniques from linear
scattering theory, such as the long time parametrix construction of Isozaki-Kitada
type \cite{IsKi}, propagation estimates, and local decay estimates.  
  
\begin{theorem}[{Global-in-time Strichartz Estimates} \cite{Mizutani}] For any admissible pairs $(q,r)$ and $(\tilde{q}^{\prime},\tilde{r}^{\prime})$ , we have 
\begin{equation}\label{strichartz1}
	\|e^{-itH}f \|_{L_t^{q}(\R:L^r_x)}\lesssim \|f\|_{L^2}
\end{equation}
and 
\begin{equation}\label{strichartz2}
	\left\|\int_0^t e^{i(t-s) H} F(s) d s\right\|_{L_t^q\left(\mathbb{R}, L_x^r\right)} \leq C\|F\|_{L_t^{\tilde{q}^{\prime}}(\R;L^{\tilde{r}^{\prime}}_x)}
\end{equation}	
	
\end{theorem}

\begin{remark}
This is the most crucial ingredient in our proof. We cannot deal with attractive Coulomb potential since the global-in-time Stricharz estimate does not hold in the attractive case.
\end{remark}

\section{Construction of the asymptotic profile $u_p$}
The construction of the asymptotic profile $u_p$ is essentially the same as the one in \cite{KaMi}. We outline below for the reader's convenience.

Given a scattering datum $u_{+}$, we define the asymptotic profile $u_p$ as follows. First, since we assume $0\notin \supp\asy$, there exists a constant $c_0>0$, such that 
$\supp\asy\subset\{|\xi|\geq c_0  \}$. We fix a cut-off radial function $\chi\in C_0^{\infty}(\R^3)$ such that $0\leq \chi\leq 1$, $\chi(x)=1$ for $|x|\leq c_0/4 $ and $\chi(x)=0$ for $|x|\geq {c_0}/2$. And we define a time-dependent potential $V_{T_1}(t,x)$ by
\begin{equation}
	V_{T_1}(t,x)=\frac{1}{|x|}\left\{1-\chi\left(\frac{2x}{t+T_1}\right) \right\},
\end{equation} 
 where $T_1>1$. Note that $V_{T_1}(t,x)=\frac{1}{|x|}$ for $t\geq 0$ and $|x|\geq c_0(t+T_1)/4$. Moreover,
 \begin{equation}
 	|\p_x^{\alpha}V_{T_1}(t,x)|\leq C_{\alpha}\langle t\rangle^{-1-|\alpha|}, \quad t\in\R,\quad x\in\R^3
 \end{equation}
 Using this decaying condition, in \cite{KaMi} it is proved that for sufficiently large $T_1>1$, there exists a solution $\Psi\in C^{\infty}([1,\infty)\times\R^3;\R)$ to the following Hamilton-Jacobi equation
 \begin{equation}\label{HamiltonJacobieq}
 	-\partial_t\Psi(t,x)=\frac{1}{2}|\nabla\Psi(t,x)|^2+V_{T_1}(t,x).
 \end{equation}

\begin{lemma}[see \cite{KaMi}] For sufficiently large $T_1\geq 1$, there exists a solution $\Psi\in C^{\infty}([1,\infty)\times\R^3;\R)$ to \eqref{HamiltonJacobieq} such that, for all $t\geq 1$, 
\begin{equation}
\left\|\nabla\Psi(t,x)-\frac{x}{t} \right\|	_{L_x^{\infty}}\lesssim \langle t\rangle^{-1},\qquad \left\|\Delta\Psi(t,x)-\frac{3}{t} \right\|_{L^{\infty}_x}\lesssim \langle t\rangle^{-2}
\end{equation}
	
\end{lemma}

Then the asymptotic profile $u_p$ is defined by 
\begin{equation}\label{asymprfile}
u_p(t,x)=(it)^{-\frac{3}{2}}e^{i\Psi(t,x)}e^{-i \{(-\Delta)^{-1}|\asy|^2)(\frac{x}{t})\log t\}}\asy(\frac{x}{t})
 \end{equation}
which can also be expressed by
\begin{equation}\label{equationup}
	u_p(t,x)=\MP(t)D(t)W(t,x)
\end{equation}

\begin{remark}
	The asymptotic profile for the critical power type nonlinearity is (see \cite{KaMi}):
\begin{equation}
	\tilde{u}_p(t,x)=(it)^{-\frac{3}{2}}e^{i\Psi(t,x)}e^{-i |\asy(\frac{x}{t})|^{\frac{2}{3}}\log t\}}\asy(\frac{x}{t})
\end{equation}
	
\end{remark}
It can also be written as:
\begin{equation}
	\tilde{u}_p(t,x)=\MP(t)D(t)\widetilde{W}(t,x)
\end{equation}

If we replace the nolinearity $F(u)$ by $\widetilde{F}(u)=|u|^{2/3}u$,  $\widetilde{W}(t,x)$  satisfies similar properties as in Lemma \ref{Wequation}. Thus, most computations in \cite{KaMi} remain the same due to Lemma \ref{Wequation}.
\section{Proof of the main theorem}
Following \cite{KaMi}, we will divide the proof of Theorem \ref{MainTheorem} into 3 steps:

\subsection{Reformulation into a fixed point argument} If we assume that $u$ is a smooth solution to \eqref{HCoulomb} satisfying $\|u-u_p \|_{L^2}\rightarrow 0$ as $t\rightarrow +\infty$. Then we can get the Duhamel formula:

\begin{equation}\label{integralequation}
	u(t)=u_p(t)+i\int_{t}^{\infty} e^{-i(t-s)H}\{F(u(s))-F(u_p(s))
    -(i\partial_s-H)u_p(s)+F(u_p(s))\} ds
\end{equation}

Although \eqref{integralequation} is different from the usual integral equation for \eqref{HCoulomb}, Lemma \ref{equivlemma} shows that a solution to \eqref{integralequation} is, in fact, a solution to the standard integral equation associated with \eqref{HCoulomb}.

We denote
$$K[u]:=i\int_t^{\infty} e^{-i(t-s)H}\{F(u(s))-F(u_p(s))\} ds$$
And the error term:

$$\mathcal{E}(t):=i\int_t^{\infty} e^{-i(t-s)H}\{-(i\partial_s-H)u_p(s)+F(u_p(s))\} ds$$
So, the integral equation \eqref{integralequation} becomes:
\begin{equation}\label{UKE}
	u(t)=u_p(t)+\mathcal{K}[u](t)+\mathcal{E}(t).
\end{equation}
Given $T,R>0$, we define the energy space $X(q, r, T,R)$ by 
\[X(q,r, T,R):=\{u\in C([T,\infty); L^2(\R^3))| \| u-u_p\|_{X_{T}}\leq R \}   \]
$$ \|u\|_{X}:=\sup_{t\geq T}t^b\|u(t)\|_{L^2} +\sup_{t\geq T} t^b\|u\|_{L^q([t,\infty);L^r)}
$$
where $(q,r)$ is admissible.

\begin{lemma}\label{equivlemma}
	Suppose $u\in C([t_0,\infty);L^2(\R^3))\cap L^q([t_0,\infty);L^r(\R^3))$ for any admissible pair $(q,r)$ and $u$ is a solution to \eqref{integralequation}   with some $t_0\in \R$. Then 
\begin{equation}\label{IntegralEquation}
	\begin{split}
		u(t)=e^{-i(t-t_0)H}u(t_0)-i\int_{t_0}^t e^{-i(t-s)H}F(u(s)) ds, \quad t\geq t_0
	\end{split}
\end{equation}		
\end{lemma}
\begin{proof}
It suffices to check the right-hand side of \eqref{IntegralEquation} and \eqref{integralequation} coincide. Changing variables $t\rightarrow t+t_0$, it is equivalent to that the following holds for $t\geq 0$:
\begin{equation}\label{UPFE}
	\begin{split}
		0= & u_{\mathrm{p}}\left(t+t_0\right)-i \int_{t+t_0}^{\infty} e^{-i\left(t+t_0-s\right) H} F\left(u_{\mathrm{p}}(s)\right) d s+\mathcal{E}\left(t+t_0\right) \\
& -e^{-i t H} u\left(t_0\right)+i \int_{t_0}^{\infty} e^{-i\left(t+t_0-s\right) H} F(u(s)) d s
	\end{split}
\end{equation}	
By \eqref{UKE} and \eqref{DecompostionMD} and \eqref{DefC}, the sum of the first three terms of the RHS of the above equation is:
\begin{equation}
\begin{aligned}
& u_{\mathrm{p}}\left(t+t_0\right)-i \int_{t+t_0}^{\infty} e^{-i\left(t+t_0-s\right) H} F\left(u_{\mathrm{p}}(s)\right) d s+\mathcal{E}\left(t+t_0\right) \\
& =\mathcal{U}_3\left(t+t_0\right) W\left(t+t_0\right)-i e^{-i\left(t+t_0\right) H} \int_{t+t_0}^{\infty} e^{i s H}\left\{\mathcal{U}_3(s) \frac{F(W(s))}{s}+\mathcal{C}(s) W(s)\right\} d s \\
& =\mathcal{U}_3\left(t+t_0\right) W\left(t+t_0\right)-i e^{-i\left(t+t_0\right) H} \int_{t+t_0}^{\infty} i \partial_s\left(e^{i s H} \mathcal{U}_3(s) W(s)\right) d s
\end{aligned}
\end{equation}
By using \eqref{UKE} with $t=t_0$, we find that the sum of the last two terms of RHS of \eqref{UPFE} is equal to: 
\begin{equation*}
\begin{split}
& -e^{-i t H} u\left(t_0\right)+i \int_{t_0}^{\infty} e^{-i\left(t+t_0-s\right) H} F(u(s)) d s \\
& =-e^{-i t H} u_{\mathrm{p}}\left(t_0\right)-i \int_{t_0}^{\infty} e^{-i\left(t+t_0-s\right) H} F\left(u_{\mathrm{p}}(s)\right) d s-e^{-i t H} \mathcal{E}\left(t_0\right) \\
& =-e^{-i t H} \mathcal{U}_3\left(t_0\right) W\left(t_0\right)-i e^{-i\left(t+t_0\right) H} \int_{t_0}^{\infty} i \partial_s\left(e^{i s H} \mathcal{U}_3(s) W(s)\right) d s
\end{split}
\end{equation*}
Then, the result follows.
\end{proof}

\subsection{Estimates of $\mathcal{K}[u]$} We prove that \eqref{integralequation} is a contraction map from $X$ to $X$ for some chosen $T, R$.

\begin{lemma}[Contraction Estimates]\label{Contractionestimates}
Suppose $b>\frac{1}{4}$ and $u,u_1,u_2\in X$
\begin{align*}
	&\|\mathcal{K}[u]\|_{X}\lesssim R\| \li|\asy|^2\|_{L^{\infty}}+T^{\frac{1}{2}-2b}R^3 \\
	&\|\mathcal{K}[u_1]-\mathcal{K}[u_2] \|_X\lesssim \|u_1-u_2\|_X(\|\li|\asy|^2\|_{L^{\infty}}+T^{\frac{1}{2}-2b}R^2)
\end{align*}	
\end{lemma}
\begin{proof}
We only show the first line and the second line follows similarly. 
Since 
$$ |F(u)-F(v)|\lesssim (\li |u|^2+\li |v|^2)|u-v|
$$	
We can decompose 
$$ F(u(s))-F(u_p(s))=F^{(1)}(u)+F^{(2)}(u)
$$
Where
\begin{align*}
	F^{(1)}(u)=1_{\{|u_p|>|u-u_p|\}}(F(u)-F(u_p)) \\
	F^{(2)}(u)=1_{\{|u_p|\leq|u-u_p|\}}(F(u)-F(u_p))
\end{align*}
So that
\begin{align*}
	&|F^{(1)}(u)|\lesssim |\li|u_p|^2|(|u-u_p|) \\
	&|F^{(2)}(u)|\lesssim |\li|u-u_p|^2|(|u-u_p|)
\end{align*}

Then the Strichartz estimate \eqref{strichartz2} gives us for any admissible pair $(q,r)$
$$ \|\mathcal{K}[u]\|_{L_t^qL_x^r}\lesssim \|\li|u_p|^2(|u-u_p)|\|_{L^1_tL^2_x}+\|\li |u-u_p|^2(|u-u_p|)\|_{L_t^\frac{4}{3}L_x^\frac{3}{2}}
$$
By scaling, we have 
\[\li|u_p|^2=t^{-1}\li |\asy|^2\]
So, the first term can be bounded easily:
\begin{align*}
	 \|\li|u_p|^2(|u-u_p)|\|_{L^1_tL^2_x} &\lesssim \int _t^{\infty}\|\li |\asy|^2\|_{L^\infty}s ^{-1} \|u-u_p\|_{L^2} ds \\
	&\lesssim t^{-b}\|\li |\asy|^2\|_{L^\infty} \|u-u_p\|_{X} \\
\end{align*}

By H\"{o}lder's inequality and Hardy-Littlewood-Sobolev inequality, we can bound the second term:
\begin{align*}
	\|\li |u-u_p|^2(|u-u_p|)\|_{L_t^\frac{4}{3}L_x^\frac{3}{2}}
 &\lesssim\| \|\li |u-u_p|^2\|_{L^6} \|u-u_p\|_{L^2} \|_{L^{\frac{4}{3}}_t} \\
	&\lesssim \| \|(u-u_p)^2\|_{L^{\frac{6}{5}}} \|u-u_p\|_{L^2}\|_{L^\frac{4}{3}} \\
	&\lesssim \| \|u-u_p\|_{L^2} \|u-u_p \|_{L^3} \|u-u_p\|_{L^2} \|_{L^\frac{4}{3}_t} \\
	&\lesssim R^2 \| s^{-2b}\|_{L^2_t} \| u-u_p\|_{L^4_tL^3_x} \lesssim R^2 t^{\frac{1}{2}-3b}\|u-u_p\|_X
\end{align*}

\end{proof}

\subsection{Control the error term $\mathcal{E}(t)$} Our goal is to show that  
for sufficiently large $T$, $\mathcal{E}(t)$ is so small that we map $X_{T,R}$ to itself.
Recall that 
\begin{equation*}
	\mathcal{E}(t):=i\int_t^{\infty} e^{-i(t-s)H}\{-(i\partial_s-H)u_p(s)+F(u_p(s))\} ds
\end{equation*}
Using \eqref{EquationW}, \eqref{EquationUp} and \eqref{equationup} , we can get 
\begin{align}
		-(i\partial_s-H)u_p(s)+F(u_p(s))
	&=e^{-isH}(-i\partial_se^{isH}M_{\Psi}D(t)W)+M_{\Psi}D(t)i\partial_s W   \label{U1U2} \\
	&=-e^{-isH}[i\partial_s,e^{isH}M_{\Psi}D(t)]W \label{U3}
\end{align}
The idea is to decompose $M_{\Psi}(t)D(t)$ into the {low-velocity parts} and {high-velocity parts}. Since 
$$ M_{\Psi}(t)D(t)=M_{\Psi}(t)D(t)\mathcal{F}\{1-M(t)\}\mathcal{F}^{-1}+M_{\Psi}D\mathcal{F}M\F^{-1}
$$
$$ M_{\Psi}D\F M\F^{-1}=\chi(\frac{x}{t})M_{\Psi}D\F M\F^{-1}+(1-\chi(\frac{x}{t}))M_{\Psi}D\F M\F^{-1}
$$
we set 	
\begin{align*}
 &\mathcal{U}_1(t):=M_{\Psi}(t)D(t)\mathcal{F}\{1-M(t)\}\mathcal{F}^{-1} \\
&\mathcal{U}_2(t):=\chi(\frac{x}{t})M_{\Psi}D\F M\F^{-1} \\
&\mathcal{U}_3(t):=(1-\chi(\frac{x}{t}))M_{\Psi}D\F M\F^{-1}
\end{align*}
where $\mathcal{U}_1$ and $\mathcal{U}_2$ are the low-velocity parts, which can be seen from the cut-off function $\chi_t(x):=\chi(x/t)$. Then we have 
\begin{equation}\label{DecompostionMD}
	 M_{\Psi}(t)D(t)=\mathcal{U}_1(t)+\mathcal{U}_2(t)+\mathcal{U}_3(t) 
\end{equation}
The operators are well understood in energy space $X$:

\begin{lemma} \label{ControlU1U2}
	Let $2b<\delta\leq 2, T>2T_1, \varepsilon>0$, we have 
	\begin{align*}
		&\|\U_1(t)f\|_{X_T}\lesssim T^{b-\frac{\delta}{2}+\varepsilon}\| f\|_{H^{\delta}} \\
		&\|\U_2(t)f\|_{X_T}\lesssim T^{b-\frac{\delta}{2}+\varepsilon} \|f \|_{H^{\delta}} \quad \mathrm{if} \quad \mathrm{supp}(f)\subset\{|x|\geq c_0 \}
	\end{align*}	
\end{lemma}
\begin{proof}
 Using Hausdorff-Young inequality, H\"{o}lder's inequality and the following:
\[|e^{\frac{i|x|^2}{2t}}-1|\lesssim \left(\frac{|x|^2}{2t}\right)^{\delta/2} 	\]
We have 
\begin{align*}
	\|\U_1(t)f\|_{L^r}&\leq \|D(t)\F\{1-M(t)\}\F^{-1}f\|_{L^r} \\
	&\lesssim |t|^{-\frac{2}{q}}\|\{1-M(t)\}\F^{-1}f \|_{L^{r^{\prime}}} \\
	&\lesssim |t|^{-\frac{2}{q}-\frac{\delta}{2}+\frac{1}{q}+\frac{\varepsilon}{2}}\||x|^{\delta-\frac{2}{q}-\varepsilon} \F^{-1}f\|_{L^{r^{\prime}}} \\
	&\lesssim |t|^{-\frac{1}{q}-\frac{\delta}{2}+\frac{\varepsilon}{2}} \||x|^{\delta-\frac{2}{q}-\varepsilon} \langle x \rangle^{\frac{2}{q}+\varepsilon} \F^{-1}f\|_{L^2}\|\langle x\rangle^{-\frac{2}{q}-\varepsilon}\|_{L^{\frac{3q}{2}}} \\
	&\lesssim |t|^{-\frac{1}{q}-\frac{\delta}{2}+\frac{\varepsilon}{2}} \|f\|_{H^{\delta}}
		\end{align*}
Taking the $L^q([t,\infty)$-norm, multiplying by $|t|^b$ and taking the supremum over $t\geq T$, we can get $X$-norm control. To show the bounds for $\U_2$, we write
\begin{align*}
	\U_2(t)f=\chi(x/t)M_{\Psi}Df+\chi(x/t)M_{\Psi}D\F\{M(t)-1\}\F^{-1}f
\end{align*}
Since we assume that $\supp(f)\subset\{|x|\geq c_0 \}$, the first term of the above vanishes if $t\geq T\geq 2T_1$.  The second term follows from the same argument as for $\mathcal{U}_1$.

\end{proof}
By  Lemma \ref{NonlinearEstimates} and Lemma \ref{ControlU1U2} , we have
\begin{align*}
\|	\mathcal{U}_{j}(t)W(t) \|_{X_T} \lesssim T^{b-\frac{1}{2}+\varepsilon}\| W(t) \|_{H^{1}} \lesssim T^{b-\frac{1}{2}+\varepsilon^{\prime}} \|\asy \|_{H^1}(1+\|\asy \|_{H^1}^2)
\end{align*}
for $j=1,2$. It follows that 
\begin{align*}
	\|	\mathcal{U}_{j}(t)W(t) \|_{L^2}\rightarrow 0 \quad \mathrm{as}\,\,t\rightarrow +\infty
\end{align*}
Therefore we can integrate \eqref{U1U2} to get 
\begin{align*}
	\mathcal{E}_{1}(t):=&\int_t^{\infty} e^{-i(t-s)H}e^{-isH}(-i\partial_se^{isH}(\U_1+\U_2)W)+e^{-i(t-s)H}(\U_1+\U_2)(i\partial_s W) \\
	&= -(\U_1(t)+\U_2(t))W+\int_t^{\infty} e^{-i(t-s)H}(\U_1+\U_2)\frac{F(W)}{s} ds
\end{align*}

\begin{lemma} For sufficiently large $T>1$, and any $\varepsilon>0$,
$$	\| \mathcal{E}_{1}\|_{X_T}\lesssim T^{b-\frac{1}{2}+\varepsilon^{\prime}}(\|\asy\|_{H^1}(1+\|\asy\|_{H^1}^2+\|\asy\|_{H^1}^4) $$
\end{lemma}
\begin{proof}
It suffices to consider the second term of $\mathcal{E}_{1}(t)$. Using Stricharz estimates with $(q^{\prime},r^{\prime})=(1,2)$, we have 
\begin{align*}
	\left\|\int_t^{\infty} e^{-i(t-s)H}\U_j(s)\frac{F(W(s))}{s} ds\right\|_X&\lesssim \sup_{t\geq T} t^b\int_t^{\infty} s^{-1} \|\U_j(s) F(W(s))\|_{L^2} ds \\
\end{align*}
By  Lemma \ref{NonlinearEstimates} and  Lemma \ref{ControlU1U2}, we have  
\begin{align*}
	\sup_{t\geq T} t^b\int_t^{\infty} s^{-1} \|\U_j(s) F(W(s))\|_{L^2} ds &\lesssim T^{b-\frac{1}{2}+\varepsilon}\| F(W(s)) \|_{H^1} \\
	&\lesssim T^{b-\frac{1}{2}+\varepsilon^{\prime}}(\|\asy\|_{H^1}(1+\|\asy\|_{H^1}^2+\|\asy\|_{H^1}^4)
)
\end{align*}
\end{proof}
It remains to consider the high-velocity part:
\[ \mathcal{E}_2(t):= -i \int_t^{\infty} e^{-i(t-s) H} \mathcal{C}(s) W(s) d s\] 
where $\mathcal{C}(s)$ denotes an extension of the commutator
\begin{equation}\label{DefC}
e^{-i s H}\left[i \partial_s, e^{i s H} \mathcal{U}_3(s)\right]=e^{-i s H} i \partial_s e^{i s H} \mathcal{U}_3(s)-i \mathcal{U}_3(s) \partial_s .
\end{equation}

The following lemma gives us the desired result:
\begin{lemma}[See \cite{KaMi}]\label{ControlU3}
 For sufficiently large $T$ and any $\varepsilon>0, e^{-i t H}\left[i \partial_t, e^{i t H} \mathcal{U}_3(t)\right]$ defined on $C_0^{\infty}((T, \infty) \times$ $\left.\mathbb{R}^3\right)$ extends to a bounded operator  $\mathcal{C}$ from $L^{\infty}([T, \infty) ; H^1)$ to $ L^1([T, \infty) ; L^2))$ which satisfies
$$
\|\mathcal{C}(t) f\|_{L^1\left([T, \infty) ; L^2\right)} \lesssim T^{-1+\varepsilon}\|f\|_{L^{\infty}\left([T, \infty) ; H^1\right)} .
$$

Moreover, 
$$
\left\|\mathcal{E}_2\right\|_{X_T} \lesssim T^{b-1+\varepsilon}  (\|\asy\|_{H^1}(1+\|\asy\|_{H^1}^2+\|\asy\|_{H^1}^4)
$$
\end{lemma}

\begin{proof}
	 This is the same as Lemma 4.4 in \cite{KaMi} since all the operators are the same.	 
\end{proof}
Now, we are ready to prove Theorem \ref{MainTheorem}.
\begin{proof}[Proof of Theorem 1.1]
Using Lemma \ref{Contractionestimates}, Lemma \ref{ControlU1U2} and Lemma \ref{ControlU3}, for any fix $R>0$, we can choose sufficiently large $T$ and  $\| \li|\asy|^2\|_{L^{\infty}}$ sufficiently small such that \eqref{integralequation} is a contraction map from $X_{T}$ to itself. Thus we have a solution $u(t,x)$ defined on $(T,\infty)$ and we extend the solution to the whole time interval by global existence result (see \cite{ChGlassey}).
\end{proof}

\bibliographystyle{amsplain}
\providecommand{\bysame}{\leavevmode\hbox to3em{\hrulefill}\thinspace}
\providecommand{\MR}{\relax\ifhmode\unskip\space\fi MR }
\providecommand{\MRhref}[2]{%
  \href{http://www.ams.org/mathscinet-getitem?mr=#1}{#2}
}
\providecommand{\href}[2]{#2}

\end{document}